\documentclass[12pt,a4paper,reqno]{amsproc}

\usepackage{amsthm}
\usepackage{amsfonts}
\usepackage{amsmath}

{ \theoremstyle{plain}
\newtheorem{thm}{Theorem}
\newtheorem{cor}{Corollary}
}

{ \theoremstyle{definition}
}

\theoremstyle{plain}

{\theoremstyle{remark}
\newtheorem{rem}{Remark}}

\theoremstyle{remark}

\def\Com#1{\mathbb{C}^{#1}}
\def\R#1{\mathbb{R}^{#1}}

\def\ov#1{\overline{#1}}

\DeclareMathOperator{\Res}{\mathcal{R}}

\begin{document}

\title{Ullemar's formula for the moment map, II}
\thanks{The paper was supported by grant RFBR
no.~03-01-00304 and grant of the Royal Swedish Academy of Sciences}

\author{Tkachev V.G.}

\email{tkatchev@math.kth.se}


\address{
Matematiska Institutionen, Kungliga Tekniska H\"ogskolan,
Lindstedtsv\"agen 25, 10044 Stockholm, Sweden }

\begin{abstract}
We prove the complex analogue of Ullemar's formula for the Jacobian
of the complex moment mapping. This formula was previously
established in the real case.
\end{abstract}

\maketitle

\section{Introduction}

Consider the `moment map'
$$
\mu:\;\Omega\to (\mu_0,\mu_1,\mu_2,\ldots),
$$
where
\begin{equation}\label{def:moment2}
\mu_k=\frac{i}{2\pi }\iint\limits_{\Omega}\zeta^k \,d\zeta\wedge
d\bar{\zeta},
\end{equation}
and $\Omega$ is a bounded domain in $\Com{}$. If $\Omega$ is a
simply-connected domain, we can uniformize it as the image
$\Omega=\phi(\mathbb{D})$, where $\phi$ is a unique function which
is holomorphic in the unit disk ${\mathbb{D}}$ and normalized by the
following conditions
\begin{equation}\label{normal}
\phi(0)=0, \qquad \phi'(0)>0.
\end{equation}
Then (\ref{def:moment2}) takes the form
\begin{equation}\label{def:moment1}
\mu_k(\phi)=\frac{i}{2\pi }\iint\limits_{\mathbb{D}}\phi^k(z)
|\phi'(z)|^2 \,dz\wedge d\bar{z}.
\end{equation}
In general, when the function $\phi$ is not globally univalent in
$\mathbb{D}$ and satisfies (\ref{normal}), we use the previous
formula as the definition for the complex (or analytic) moments of
 $\phi$ \cite{BG84}; then $\Omega$ is regarded as a Riemannian surface
over $\Com{}$.

This notion appears in several problems of complex analysis and its
applications. In particular, the sequence $(\mu_k)_{k\geq 1}$
constitutes an infinite family of invariants of the Hele-Shaw
problem \cite{R72} of the cell $\Omega$ and can be used as a canonic
coordinate system in the corresponding Laplacian growth model
\cite{Zab} (see also \cite{P} for the functional analysis
interpretation).

In what follows, we consider the special case when  the moment map
$\mu$ is restricted to the set $\mathfrak{S}_n\subset\Com{}[z]$ of
polynomials of degree $n$ normalized by (\ref{normal}). It is easy
to verify that $\mu_k=0$ for all $k\geq n$ (see also (\ref{zero})
below), so only the first $n$ moments are of interest for $P\in
\mathfrak{S}_n$. We consider the induced  finite dimensional map
\begin{equation}\label{def-momentmap}
\mu:\; P=a_1z+a_2z^2+\ldots +a_nz^n \to
(\mu_0(P),\ldots,\mu_{n-1}(P)).
\end{equation}
Since $\mu_0(P)>0$,
$$
\mu:\; \mathfrak{S}_n\to \R{+}\times\Com{n-1}.
$$
Notice that
$$
\dim_{\R{}} \R{+}\times\Com{n-1}=\dim_{\R{}}\mathfrak{S}_n=2n-1,
$$
where $\mathfrak{S}_n$ is regarded as an open subset of a linear
space. In \cite{BG84} Gustafsson proved that the Fr\'{e}chet
derivative $d\mu$ is non-singular at any $P$ which is a locally
univalent polynomial.

On the other hand, notice that the subset  $\mathfrak{S}^{\R{}}_n$
which consists of the polynomials with real coefficients is an
invariant set of $\mu$  in the sense that all $\mu_k(P)$ are real
(cf. (\ref{summa}) below). Hence, in a similar way  the moment map
induces the map $\mu_{\R{}}:\mathfrak{S}^{\R{}}_n\to\R{n}$.

In \cite{Ul}, C.~Ullemar conjectured the following formula for the
Jacobian of $\mu_{\R{}}$:
\begin{equation}\label{J}
J_{\R{}}(P):=\frac{\partial (\mu_0,\ldots,\mu_{n-1})}{\partial
({a}_{1},\ldots,a_{n})}
(P)=2{}^{-\frac{n(n-3)}{2}}a_1{}^{\frac{n(n-1)}{2}}P'(1)P'(-1)\Delta
({P'}^*(z)),
\end{equation}
where $\Delta$ stands for the principal Hurwitz determinant
\cite[\S15.715]{GR}, and $Q^*$ denotes the mirror conjugate image of
polynomial $Q$, i.e.
\begin{equation}\label{eee}
Q^{*}(z):= z^{m}\bar{Q}(1/z)=\bar{q}_{m}+\bar{q}_{m-1}z+\ldots
+\bar{q}_0z^m, \qquad m=\deg Q,
\end{equation}
where $\bar{Q}(z)=\overline{Q(\bar{z})}$ is the conjugate
polynomial.

The above expression for the Jacobian was recently proved in
\cite{KT} as a consequence of the following identity
\begin{equation}\label{-A31}
J_{\R{}}^2(P)=4(-1)^{n-1}a_1^{n(n-1)}\Res (P',{P'}^*)\cdot
P'(-1)P'(1),
\end{equation}
where $\Res (\cdot,\cdot)$ denotes the resultant of the
corresponding polynomials.

In this paper we generalize formula (\ref{-A31}) for polynomials
with arbitrary complex coefficients.

\begin{thm}\label{th1}
The Jacobian of the moment map is expressed as follows
\begin{equation}\label{main-formula}
J_{\Com{}}(P):=\frac{\partial (\bar{\mu}_{n-1},\ldots,\bar{\mu}_{1},
\mu_0, \mu_1,\ldots,\mu_{n-1})}{\partial
(\bar{a}_{n},\ldots,\bar{a}_{2}, a_1,
a_2,\ldots,a_{n})}=2a_1^{n^2-n+1}\Res(P',P'^*).
\end{equation}
\end{thm}

A well known theorem of Sylvester (see Section~\ref{sub-res}) allows
us to compute the above resultant as the determinant of a matrix of
size $2n-2$, whose entries are $0$ or a coefficient of either $P'$
or $P'^*$. In particular, the resultant is homogeneous in the
coefficients of $P'$ and $P'^*$ separately, with respective degree
$n-1$.

On the other hand, geometrically, the hypersurface
$$
\{(a,\bar{a}):\quad \Res(P',P'^*)=0\}
$$
i.e. the critical set of the Jacobian, is the projection of the
incidence variety
$$
\{(a,\bar{a},z):\quad
\sum_{k=1}^nka_kz^{k-1}=\sum_{k=0}^{n-1}(n-k)\bar{a}_{n-k}z^k\}
$$
that is to say, the set of $(a,\bar{a})$ which appear above for some
$z$.

The following assertion is a direct consequence  of the definition
(\ref{resultant-def}) of the resultant and formula (\ref{eee})
above, and it characterizes the set of critical points of $d\mu$.

\begin{cor}
The moment map is degenerate at $P$ if and only if the derivative
$P'$ has two roots $\alpha_i$ and $\alpha_j$ such that
$\alpha_i\bar{\alpha}_j=1$ (the case $i=j$ is permitted).
\end{cor}

Note that for a locally univalent polynomial in the closed unit disk
we have $|\alpha_j|<1$ for all the roots of its derivative. Hence,
we obtain another proof of the above result due to Gustafsson
\cite{BG84}.

\begin{cor}
The moment map is locally injective on the set of all locally
univalent polynomials in the closed unit disk.
\end{cor}


\section{Preliminaries}\label{s1}

\subsection{Complex moments}

Using the Stokes formula, we obtain
\begin{equation}\label{stockes}
\mu_k=\frac{i}{2\pi (k+1)}\int\limits_{\partial\Omega}\zeta^{k+1} \,
d\bar{\zeta}=\frac{1}{2\pi
i}\int\limits_{\partial\Omega}\zeta^{k}\bar{\zeta} \, d\zeta,
\end{equation}
which implies
\begin{equation*}\label{m-rep0}
\mu_k(\phi)=\frac{i}{2\pi (k+1)}\int\limits_{
\mathbb{T}}\phi^{k+1}(z) \bar{\phi}'(\bar{z})\,
d\bar{z}=\frac{1}{2\pi i}\int\limits_{\mathbb{T}}\phi^{k}(z)
\bar{\phi}(\bar{z}) \phi'(z)\, dz,
\end{equation*}
where $\mathbb{T}=\partial \mathbb{D}$ is the unit circle. Hence,
using the identity $\bar{z}=1/z$ which holds everywhere in
$\mathbb{T}$, we get
\begin{equation}\label{m-rep1}
\mu_k(\phi)=\frac{1}{2\pi i(k+1)}\int\limits_{
\mathbb{T}}\phi^{k+1}(z) \bar{\phi}'(1/z)\,
\frac{dz}{z^2}=\frac{1}{2\pi i}\int\limits_{\mathbb{T}}\phi^{k}(z)
\bar{\phi}(1/z) \phi'(z)\, dz.
\end{equation}

Given a function which is analytic in a neighborhood of
$\mathbb{T}$, let us denote by $\lambda_s(f)$ the $s$th Laurent
coefficient of $f$, i.e.
$$
f(z)=\sum_{s=-\infty}^{\infty}\lambda_s(f)z^s,
$$
hence
\begin{equation}\label{lambda}
\mu_k(\phi)=\frac{1}{k+1}\lambda_{1}(\phi^{k+1}(z)\bar{\phi}'(1/z))
=\lambda_{-1}(\phi^{k}(z) \phi'(z)\bar{\phi}(1/z) )
\end{equation}

Now, let $P$ be an arbitrary polynomial in $\frak{S}_n$. Then
$\bar{P}'(1/z)=P'^*(z)z^{1-n}$ and $\bar{P}(1/z)=P^*(z)z^{-n}$,
which by virtue of (\ref{lambda}) yields
\begin{equation}\label{lambda1}
\mu_k(P)=\frac{1}{k+1}\lambda_{n}(P^{k+1}P'^*)
=\lambda_{n-1}(P'P^{k}P^* )
\end{equation}

It follows from the first identity in (\ref{lambda1}) and $P(0)=0$
that
\begin{equation}\label{zero}
\mu_k(P)=0, \qquad k\geq n
\end{equation}
On the other hand, the second identity in (\ref{lambda1}) yields the
so-called \textit{Richardson formula}
\begin{equation}
\mu_k(P)=\sum s_1 a_{s_1}\cdots a_{s_{k+1}}\bar{ a}_{s_1+\ldots
+s_{k+1}}, \label{summa}
\end{equation}
where the sum is taken over all possible sets of indices $s_1$,
$\ldots $, $s_k\geq 1$. It is assumed that $a_j=0$ for $j\geq n+1$.
These formulae are easy to use for straightforward manipulations
with the complex moments and it follows also that $\mu_k(P)$ is a
\textit{polynomial} mapping.

It is convenient to identify $\frak{S}_n$ with the corresponding
coefficient subset in $\R{+}\times\Com{n-1}$ in a standard way:
$$
a\sim P:=a_1z+a_2z^2+\ldots+a_nz^n.
$$
Since,
$$
\mu_0(P)=\sum_{s=1}^{n}s|a_s|^2>0, \qquad
\mu_{n-1}(P)=na_{1}^{n}\bar{a}_n\ne 0,
$$
the moment  map (\ref{def-momentmap}) is well defined as an
automorphism of $\frak{S}_n$ into itself.

\subsection{Resultants}\label{sub-res}
Here we review some basic facts about the resultant; see
\cite{WanDer} for a detailed introduction.

The \textit{resultant} of two polynomials
$$
A(z)=a_m\prod_{j=1}^m(z-\alpha_j), \qquad
B(z)=b_k\prod_{j=1}^k(z-\beta_j)
$$
with respect to $z$ is the polynomial
\begin{equation}\label{resultant-def}
    \Res(A,B)=a_m^kb_k^m\prod_{i,j=1}(\alpha_i-\beta_j).
\end{equation}
The resultant vanishes iff $A$ and $B$ have a common root. It can be
evaluated as the determinant of the \textit{Sylvester matrix}, which
is the following $m+k$  by $m+k$ matrix
$$
\begin{pmatrix}
  a_0 &\quad  a_1 &\quad  \ldots &\quad   \ldots    &\quad a_m     &\quad               \\
      &\quad  a_0 &\quad  a_1    &\quad \ldots&\quad    \ldots     &\quad a_m           \\
    &\quad   &\quad   &\quad  &\quad \vdots &\quad &\quad  \\
  &\quad    &\quad      &\quad  a_{0}  &\quad  a_{1} &\quad  \ldots  &\quad \ldots &\quad a_m    \\
  b_0 &\quad  b_1 &\quad  \ldots &\quad b_k     &\quad    &\quad                \\
      &\quad  b_0 &\quad  b_1    &\quad \ldots&\quad b_k &                   \\
      &\quad   &\quad &\quad &\quad \vdots  &\quad&\quad  \\
    &\quad   &\quad      &\quad  b_{0}  &\quad  b_{1} &\quad  \ldots  &\quad\ldots&\quad b_k    \\
\end{pmatrix}
$$
in which the first $k$ rows are the coefficients of $A$, the next
$m$ rows are the coefficients of $B$, and the elements not shown are
all zero. The following are some useful elementary properties we
will use below.
\begin{equation}\label{res-prop}
\begin{split}
\Res(A,B)&=(-1)^{km}\Res(B,A),\\
\Res(A_1A_2,B)&=\Res(A_1,B)\Res(A_2,B),\\
\Res(z^n,A)&=A^n(0).
\end{split}
\end{equation}

Next, given a polynomial $A(z)$ of degree $n$, we define its mirror
conjugate image as
$$
A^{*}(z):= z^{n}\bar{A}(1/z)=\bar{a}_{n}+\bar{a}_{n-1}z+\ldots
+\bar{a}_0z^n,
$$
where $\bar{A}(z)=\overline{A(\bar{z})}$ is the conjugate
polynomial.  We have for their roots: $
\alpha^*_j=(\bar{\alpha_j})^{-1}$ and the corresponding resultant
takes the following form
\begin{equation}\label{rez_defA}
\Res(A,A^*)=\det
\begin{pmatrix}
  a_0 & a_1 & \ldots & \ldots     &a_n     &              \\
      & a_0 & a_1    &\ldots&  \ldots      &a_n           \\
   &  &  & &\vdots &  & \\
  &   &     & a_{0}  & a_{1} & \ldots  &\ldots &a_n    \\
  \bar{a}_n & \bar{a}_{n-1} & \ldots & \ldots     &\bar{a}_0     &               \\
      & \bar{a}_n & \bar{a}_{n-1}    &\ldots& \ldots       &\bar{a}_0             \\
     &  & &  &\vdots &  & \\
  &   &     & \bar{a}_{n}  & \bar{a}_{n-1} & \ldots  &\ldots &\bar{a}_0    \\
\end{pmatrix}.
\end{equation}

\begin{rem}
We wish to point out that the latter form, $\Res(A,A^*)$, is
irreducible as a polynomial of $(a,\bar{a})$ over $\Com{}$. The
proof is given in \cite[Theorem~6]{KT}.
\end{rem}

\section{Proof of the Theorem}

First, we evaluate the partial derivative of the moment map. Namely,
we have for all $k=0,\ldots,n-1$, $j=1,\ldots,n$, except for $j=1$,
$k=0$,
\begin{equation}\label{id-main}
\begin{split}
\frac{\partial \mu_k(P)}{\partial a_j}&=\lambda_{n-j}(P'^*P^{k}),\\
\frac{\partial \mu_k(P)}{\partial
\bar{a}_j}& =\lambda_{j-1}(P'P^{k}).\\
\end{split}
\end{equation}

In fact, let $j$ be an integer from $\{2,\ldots,n\}$. Then by the
first identity in (\ref{lambda1}) we have for
$$
\frac{\partial \mu_k(P)}{\partial
a_j}=\frac{1}{k+1}\lambda_{n}(P'^*\frac{\partial P^{k+1}}{\partial
a_j})=\lambda_{n}(P'^*P^{k} z^j)=\lambda_{n-j}(P'^*P^{k}).
$$
Similarly, using
$P^*=\bar{a}_n+\bar{a}_{n-1}z+\ldots+\bar{a}_2z^{n-2}+a_1z^{n-1}$
and the second identity in (\ref{lambda1}) we obtain
$$
\frac{\partial \mu_k(P)}{\partial
\bar{a}_j}=\lambda_{n-1}(P'P^{k}\frac{\partial P^*}{\partial
\bar{a}_j} )=\lambda_{n-1}(P'P^{k}z^{n-j})=\lambda_{j-1}(P'P^{k}).
$$
Finally, for $j=1$ we have by the first identity in (\ref{lambda1})
$$
\frac{\partial \mu_k(P)}{\partial
a_1}=\frac{1}{k+1}\lambda_{n}(P'^*\frac{\partial P^{k+1}}{\partial
a_1}+P^{k+1}\frac{\partial P'^*}{\partial a_1})=
\lambda_{n-k}(P'^*P^{k})+\frac{1}{k+1}\lambda_1(P^{k+1}).
$$
But $\lambda_1(P^{k+1})=0$ for $k\geq 1$, hence the desired
assertion follows.

We will make use the following notation
\begin{equation*}
\nabla f:=(\frac{\partial f}{\partial a_n},\,\frac{\partial
f}{\partial a_{n-1}},\ldots,\frac{\partial f}{\partial
a_2},\,\frac{\partial f}{\partial a_1},\,\frac{\partial f}{\partial
\bar{a}_2},\ldots, \frac{\partial f}{\partial
\bar{a}_{n-1}},\,\frac{\partial f}{\partial \bar{a}_n}),
\end{equation*}
and by
$$
q_{j-1}=ja_j
$$
we denote the coefficients of the derivative $Q:=P'$. Then
\begin{equation}\label{grad-1}
\nabla
\mu_0=(\bar{q}_{n-1},\ldots,\bar{q}_{1},2q_0,{q}_{1},\ldots,{q}_{n-1}),
\end{equation}
and for all $k=1,\ldots,n-1$ we have from (\ref{id-main})
\begin{equation}\label{grad-k}
\nabla
\mu_k=(\lambda_{0}(Q^*P^{k}),\ldots,\lambda_{n-2}(Q^*P^{k}),\lambda_{n-1}(Q^*P^{k}),
\lambda_{1}(QP^{k}),\ldots,\lambda_{n-1}(Q^*P^{k})).
\end{equation}

Let $\mathbf{Y}_0=\nabla\mu_0$ and for $k\geq 1$ write
$$
\mathbf{Y}_k:=(\lambda_{0}(Q^*z^{k}),\ldots,\lambda_{n-2}(Q^*z^{k}),\lambda_{n-1}(Q^*z^{k}),
\lambda_{1}(Qz^{k}),\ldots,\lambda_{n-1}(Q^*z^{k})).
$$
As a direct consequence of the above formula we conclude that
$$
\mathbf{Y}_k=\mathbf{0}, \quad k\geq n.
$$

Then it follows from (\ref{grad-k}) and
$$
P=z(a_1+\ldots+a_nz^{n-1})
$$
that for all $k\geq1$
$$
\nabla\mu_k=a_1^{k}\mathbf{Y}_k+\sum_{j=k+1}^{n-1}w_{k,j}\mathbf{Y}_j.
$$
Thus,
\begin{equation}\label{wedge}
\nabla\mu_0\wedge \nabla \mu_1\wedge\cdots\wedge \nabla\mu_{n-1}=
a_1^{N}\;\mathbf{Y}_0\wedge \mathbf{Y}_1\wedge\cdots\wedge
\mathbf{Y}_{n-1},
\end{equation}
where $N=(n-1)n/2$.

On the other hand, for all $k\geq1$ we have
$$
\mathbf{Y}_k:=(0,\ldots,0,\bar{q}_{n-1},\ldots,\bar{q}_k,0,\ldots,0,q_0,q_1,\ldots,q_{k-1}),
$$
where the zeroes groups contain $k$  and $k-1$ items respectively.

Now we treat the conjugate moments. We have
$\bar{\mu}_k(P)=\mu_k(\bar{P})$, whence
$$
\nabla \bar{\mu}_k=(\nabla \mu_k)^*,
$$
where by $\mathbf{X}^*$ we denote the mirror conjugate image of
vector $ \mathbf{X}=(x_1,x_2,\ldots,x_{2n-1}), $ i.e.
$$
\mathbf{X}^*=(\bar{x}_{2n-1},\ldots,\bar{x}_2,\bar{x}_1).
$$
Repeating the above argument for the conjugate expressions yields
\begin{equation}\label{wedge1}
\nabla\bar{\mu}_{n-1}\wedge \nabla \bar{\mu}_{n-2}\wedge\cdots\wedge
\nabla\bar{\mu}_{1}= a_1^{N}\;\mathbf{Y}^*_{n-1}\wedge
\mathbf{Y}^*_{n-2}\wedge\cdots\wedge \mathbf{Y}^*_{1},
\end{equation}
hence
\begin{multline}
\nabla\bar{\mu}_{n-1}\wedge\cdots\wedge \nabla\bar{\mu}_{1}\wedge
\nabla\mu_0\wedge \nabla \mu_1\wedge\cdots\wedge \nabla\mu_{n-1}=\\
=a_1^{n^2-n}\;\mathbf{Y}^*_{n-1}\wedge\cdots\wedge
\mathbf{Y}^*_{1}\wedge \mathbf{Y}_0\wedge
\mathbf{Y}_1\wedge\cdots\wedge \mathbf{Y}_{n-1}.
\end{multline}

We rewrite the latter identity in terms of determinants which gives
the following expression for the Jacobian
\begin{equation}\label{det-mu}
J_{\Com{}}(P)= \frac{\partial (\bar{\mu}_{n-1},\ldots,\bar{\mu}_{1},
\mu_0, \mu_1,\ldots,\mu_{n-1})}{\partial
(\bar{a}_{n-1},\ldots,\bar{a}_{1}, a_0,
a_1,\ldots,a_{n-1})}=a_1^{n^2-n}\det\mathbf{Y},
\end{equation}
where
$$
\mathbf{Y}=
\begin{pmatrix}
q_0 &  &  & & q_{n-1} & & &  &  \\
\bar{q}_1 & q_0 &  &  & q_{n-2} &  q_{n-1} &  &  &  \\
\vdots & \vdots & \ddots &  & \vdots & \vdots & \ddots &  &  & \\
\bar{q}_{n-2} & \bar{q}_{n-3} & \ldots & \bar{q}_{0} & q_1 & q_{2} & q_{3} & \ldots & q_{n-1} &   \\
\bar{q}_{n-1} & \bar{q}_{n-2} & \ldots & \bar{q}_{1} & 2q_0 & q_{1} & q_{2} & \ldots & q_{n-2} & q_{n-1}  \\
& \bar{q}_{n-1} & \ldots & \bar{q}_2 & \bar{q}_1 & q_0  & q_1& \ldots & q_{n-3} & q_{n-2} \\
& & \ddots & \vdots &  \vdots & & \ddots &  &  \vdots & \vdots  \\
& & & \bar{q}_{n-1} &\bar{q}_{n-2} &  &  & & q_{0} & q_1  \\
& & &  &\bar{q}_{n-1} &  &  & & & q_0  \\
\end{pmatrix},
$$
and the elements not shown are all zero.

Now, let $\mathbf{X}_{j}$ denote the $j$th column in $\mathbf{Y}$.
We have for $j=1,\ldots,n-1$
$$
\mathbf{X}_{j}=(0,\ldots,0,q_0,\bar{q}_1,\ldots,\bar{q}_{n-1},0,\ldots,0)^\top,
$$
with $j-1$ first zeroes, and for $j=n+1,\ldots,2n-1$:
$$
\mathbf{X}_{j}=(0,\ldots,0,q_{n-1},\ldots,{q}_{1},q_0,0,\ldots,0)^\top,
$$
with $j-n$ first zeroes, and
$$
\mathbf{X}_n=(q_{n-1},\ldots,q_1,2q_0,\bar{q}_1,\ldots,\bar{q}_{n-1})^\top.
$$

One can readily verify that
$$
\mathbf{X}_n-\sum_{j=1}^{n-1}\frac{q_{n-j}}{q_0}\mathbf{X}_j+
\sum_{j=n+1}^{2n-1}\frac{\ov{q}_{j-n}}{q_0}\mathbf{X}_j=
(0,\ldots,0,2q_0,2\bar{q}_1,\ldots,2\bar{q}_{n-1})^\top,
$$
which yields for the determinant
\begin{equation}\label{dett}
\det \mathbf{Y}= 2\det
\begin{pmatrix}
q_0 &  &  & &  & & &  &  \\
\bar{q}_1 & q_0 &  &  &  & \quad q_{n-1} &  &  &  \\
\vdots & \vdots & \ddots &  & & \vdots &  & \ddots  & & \\
\bar{q}_{n-1} & \bar{q}_{n-2} & \ldots & \bar{q}_{1} & q_0 & q_{1} & q_{2} & \ldots & q_{n-2} & q_{n-1}  \\
& \bar{q}_{n-1} & \ldots & \bar{q}_2 & \bar{q}_1 & q_0  & q_1& \ldots & q_{n-3} & q_{n-2} \\
& & \ddots & \vdots &  \vdots & & \ddots &  &  \vdots & \vdots  \\
& & & \bar{q}_{n-1} &\bar{q}_{n-2} &  &  & & q_{0} & q_1  \\
& & &  &\bar{q}_{n-1} &  &  & & & q_0  \\
\end{pmatrix},
\end{equation}
The latter is the transposed Sylvester matrix of $Q^*(z)$ and
$zQ(z)$, hence by (\ref{res-prop})
\begin{equation*}\label{last}
\begin{split}
\det \mathbf{Y}&=2\Res(Q^*,zQ)=2(-1)^{n(n-1)}\Res(zQ,Q^*)=\\
&=2\Res(z,Q^*)\Res(Q,Q^*)=2Q(0)\Res(Q,Q^*).
\end{split}
\end{equation*}
Thus, using our notation $Q=P'$ we arrive at
$$
J_{\Com{}}(P)=2a_1^{n^2-n+1}\Res(P',P'^*),
$$
which completes the proof.

\section*{Acknowledgements}

The author wish to thank anonymous referees for valuable comments
and suggestions.


\begin{thebibliography}{99}

\bibitem{GR}
\textsc{I.S. Gradshteyn and I.M. Ryzhik.} Tables of Integrals,
Series, and Products, 6th ed. San Diego, CA: Academic Press, 2000.

\bibitem{BG84}
\textsc{B.~Gustafsson}, On a differential equation arising in a
Hele-Shaw flow moving boundary problem. \textit{ Ark. f\"or Mat.},
\textbf{22}(1984), 251--268.

\bibitem{Zab}
\textsc{I. Krichever, M. Mineev-Weinstein, P.Wiegmann, A. Zabrodin,}
Laplacian growth and Whitham equations of soliton theory, (arXiv:
nlin.SI/0311005).

\bibitem{R72}
\textsc{S.~Richardson}, Hele-Shaw flows with a free boundary
produced by the injection of fluid into a narrow channel. J. Fluid
Mech., \textbf{56}(1972), 609--618.

\bibitem{KT}
\textsc{O.S.~ Kuznetsova, V.G.~Tkachev.} Ullemar's formula for the
Jacobian of the complex moment mapping, \textit{Complex Variables
and Applications}, \textbf{49}(2004), No~ 1, 55--72.


\bibitem{P}
\textsc{M.~Putinar}, Linear Analysis of Quadrature Domains, III.
\textit{J. Math. Anal. Appl.}, \textbf{239}(1999), 101--117.

\bibitem{Ul}
\textsc{C.~Ullemar}, Uniqueness theorem for domains satisfying
quadrature identity for analytic functions. \textit{ TRITA-MAT
1980-37, Mathematics.}, (1980) Preprint of Royal Inst. of
Technology, Stockholm.

\bibitem{WanDer}
\textsc{Van Der Warden B.L.}, Modern algebra. Vol.~1. Springer.
Berlin, 1971.



\end{thebibliography}
\end{document}